# ON STEPWISE CONTROL OF THE GENERALIZED FAMILYWISE ERROR RATE

By Wenge Guo and M. Bhaskara Rao

National Institute of Environmental Health Sciences and University of Cincinnati

A classical approach for dealing with the multiple testing problem is to restrict attention to procedures that control the familywise error rate (FWER), the probability of at least one false rejection. In many applications, one might be willing to tolerate more than one false rejection provided the number of such cases is controlled, thereby increasing the ability of the procedure to detect false null hypotheses. This suggests replacing control of the FWER by controlling the probability of $k$ or more false rejections, which is called the $k$-FWER. In this article, a unified approach is presented for deriving the $k$-FWER controlling procedures. We first generalize the well-known closure principle in the context of the FWER to the case of controlling the $k$-FWER. Then, we discuss how to derive the $k$-FWER controlling stepwise (stepdown or stepup) procedures based on marginal $p$-values using this principle. We show that, under certain conditions, generalized closed testing procedures can be reduced to stepwise procedures, and any stepwise procedure is equivalent to a generalized closed testing procedure. Finally, we generalize the well-known Hommel procedure in two directions, and show that any generalized Hommel procedure is equivalent to a generalized closed testing procedure with the same critical values.

**1. Introduction.** Consider the problem of simultaneously testing a finite number of null hypotheses $H_i$ $(i = 1, \cdots, n)$, using tests that are available for these individual hypotheses. A traditional concern dealing with this problem is to control the familywise error rate (FWER), the probability of falsely rejecting at least one true null hypothesis. However, quite often, when the number $n$ of hypotheses is large, control of FWER is so stringent that a few of the false null hypotheses are rejected. Therefore, the classical idea of controlling the FWER has been recently generalized to that of controlling the probability of $k$ or more false rejections, which is called the $k$-FWER. The rationale behind the $k$-FWER is that, often in practice, one is willing to tolerate a few false rejections, so by controlling $k$ or more false rejections the ability of a procedure to detect more false null hypotheses can potentially be improved.

A number of methods controlling the $k$-FWER have been recently suggested. Single-step and stepwise (stepdown and stepup) $k$-FWER procedures in terms of only the marginal null distributions of the $p$-values under arbitrary dependence of the $p$-values are derived in Hommel and Hoffmann (1987), Lehmann and Romano (2005), and Romano and Shaikh (2006). When the $p$-values are positively dependent, Sarkar (2008a) developed alternative single-step and stepwise $k$-FWER procedures utilizing $k$th order joint null distributions of the $p$-values. When the $p$-values are independent, Guo and Romano (2007) provided more powerful single-step and stepdown $k$-FWER procedures. In Korn et al. (2004), alternative permutation-based procedures are proposed to control the $k$-FWER approximately, which account for the dependence structure of the individual test statistics or $p$-values. Their results were generalized in Romano and Wolf (2007). In van der Laan et al. (2004), alternative procedures controlling the $k$-FWER are provided by augmenting single-step and stepwise FWER procedures. Further methods are discussed in Dudoit et al. (2004) and van der Laan et al. (2005).

In contrast to the popular false discovery rate (FDR), it is easier to derive powerful $k$-FWER controlling procedures in numerous settings. For example, suppose we are examining all pairwise comparisons in the one-way ANOVA model, in which the number of treatments is moderate. In this situation, the assumption of positive regression dependence of the underlying test statistics is not satisfied (Yekutieli, 2008), so Benjamini and Hochberg (1995)'s procedure is not applicable (Benjamini and Yekutieli, 2001), but no other FDR controlling procedure is available dealing with this problem. An alternative choice is to control the $k$-FWER, since we will be able to develop relatively easy powerful $k$-FWER controlling procedures accounting for the special dependence structure of the individual test statistics. Hence, in many applications, the $k$-FWER can be regarded as a good complement to the FWER and FDR. For further enunciation of $k$-FWER criterion, see Hommel and Hoffmann (1987), Korn et al. (2004), Lehmann and Romano (2005), and van der Laan et al. (2004). In addition, based on the similar rationale to the $k$-FWER, Sarkar (2008b) advocated the $k$-FDR using the expected ratio of $k$ or more false rejections to the total number of rejections, which is a generalization of the FDR. Several procedures controlling the $k$-FDR have also been developed in Sarkar (2008b) and Sarkar and Guo (2008a, b).

In this paper, we focus on the control of the $k$-FWER. Instead of pursuing a piece-meal approach, a unified approach is provided for the construction of the $k$-FWER controlling procedures based on marginal $p$-values. The main motivation comes from one particular paradigm of research on the FWER, where the well-known closure principle plays a fundamental role in the construction of the FWER controlling procedures. We believe that a generalization of the closure principle will play a similar key role in the construction of the $k$-FWER controlling procedures. To begin with, we generalize the closure principle, and then derive several general results on the relationship between generalized closed testing procedures and stepdown, stepup, and generalized Hommel procedures. As an application, it is then shown that the existing procedures can be directly derived following the generalized closure principle, and they are equivalent to some generalized closed testing procedures.

This paper is organized as follows. In Section 4.2, we set up the terminology. A generalization of the closure principle and several global tests are provided in Section 4.3. In Section 4.4, we discuss the relationship between generalized closed testing procedures and stepdown and stepup procedures. Several general results are obtained. In Section 4.5, we generalize the Hommel procedure and show that generalized Hommel procedures are equivalent to generalized closed testing procedures with the same critical values. In Section 4.6, we offer some concluding remarks.

**2. Basic Setting.** Consider the problem of testing simultaneously a family of $n$ null hypotheses $H_1, \ldots, H_n$. Suppose that the family satisfies the *free combination condition* of Holm (1979), that is, for any $I \subseteq \{1, \ldots, n\}$, there exists a distribution $P \in \Omega$, for which all $H_i, i \in I$ are true and all $H_i, i \notin I$ are false, where $\Omega$ is the set of all possible distributions of the data.

Suppose $V$ is the number of true null hypotheses falsely rejected. The generalized familywise error rate ($k$-FWER) is defined to be the probability of at least $k$ false rejections, where $k$ is pre-specified with $1 \leq k \leq n$. That is,

$$(2.1) \qquad k\text{-FWER} = P\{V \geq k\} \ .$$

If $k = 1$, $k$-FWER is the usual familywise error rate (FWER). When testing $H_1, \ldots, H_n$, we assume that the $p$-values $P_1, \ldots, P_n$ are available, and the $p$-values associated with the true null hypotheses satisfy

$$(2.2) \qquad P\{P_i \leq u\} \leq u \quad \text{for any } u \in (0,1).$$

Let the ordered $p$-values be denoted by $P_{(1)} \leq \cdots \leq P_{(n)}$, and the associated hypotheses by $H_{(1)}, \ldots, H_{(n)}$. Suppose $\alpha_k \leq \cdots \leq \alpha_n$ be a non-decreasing sequence of critical values.

There are two main avenues open for developing multiple testing procedures based on marginal $p$-values: stepup or stepdown. We generalize these procedures to accommodate control of the $k$-FWER. A (generalized) stepup procedure based on the critical values $\alpha_i$, which is slightly different from the usual one, is described below. If $P_{(n)} \leq \alpha_n$, then reject all null hypotheses; otherwise, reject hypotheses $H_{(1)}, \cdots, H_{(r)}$, where $r \geq k$ is the smallest

index satisfying

$$(2.3) \qquad P_{(n)} > \alpha_n, \cdots, P_{(r+1)} > \alpha_{r+1}.$$

If, for all $r \geq k$, $P_{(r)} > \alpha_r$, then reject the first $(k-1)$ most significant hypotheses.

Similarly, a (generalized) stepdown procedure, which is slightly different from the usual one, is described below. If $P_{(k)} > \alpha_k$, reject the first $(k-1)$ most significant hypotheses. Otherwise, reject hypotheses $H_{(1)}, \cdots, H_{(r)}$ where $r \geq k$ is the largest index satisfying

$$(2.4) \qquad P_{(k)} \leq \alpha_k, \cdots, P_{(r)} \leq \alpha_r.$$

Note that, if $k = 1$, the stepwise (stepup or stepdown) procedures described above are the same as the usual ones.

Evidently, from the definition of the $k$-FWER, one can always reject the $(k-1)$ most significant hypotheses without violating control of the $k$-FWER. This is the reason why we give a slightly different definitions of stepup and stepdown procedures, in which, the $(k-1)$ most significant hypotheses are automatically rejected. An alternative choice is to let $\alpha_i = \alpha_k$, $1 \leq i < k$, as in Hommel and Hoffmann (1987) and Lehmann and Romano (2005). For convenience of discussion, in the subsequent sections, all procedures including the closed testing procedures described in Section 4.3 and the generalized Hommel procedures defined in Section 4.5, are also supposed to reject automatically the $(k-1)$ most significant hypotheses.

**3. Generalized Closure Principle.** In this section, we generalize the well-known closure principle in the context of the FWER to the case of controlling the $k$-FWER. Similar to the usual closure principle, the value of the generalized closure principle is that the problem of controlling the $k$-FWER is reduced to the problem of controlling the usual probability of the Type I error of single tests of intersection hypotheses.

Let $I \subset \{1, \ldots, n\}$ with $|I| \geq k$. Let $P_i, i \in I$ be the $p$-values associated with $H_i, i \in I$ and $P_{1:I} \leq \cdots \leq P_{|I|:I}$ an ordered arrangement of these $p$-values. Let $H_{1:I}, \ldots, H_{|I|:I}$ be the corresponding null hypotheses. Let $H_I = \cap_{i \in I} H_i$. Let $\alpha_{k:I} \leq \cdots \leq \alpha_{|I|:I}$ be a given set of critical values. We present the following local test based on marginal $p$-values, which is proposed in Sarkar (2008a), for testing the single hypothesis $H_I$:

(3.1) *Reject $H_I$ if and only if $P_{j:I} \leq \alpha_{j:I}$ for at least one $j \in \{k, \ldots, |I|\}$,*

which implies that an intersection hypothesis is declared significant if and only if at least $k$ of the individual hypotheses is found significant. The local test is denoted by $(I : \alpha_{k:I}, \ldots, \alpha_{|I|:I})$ and the Type I error probability associated with the local test is given by $P(\cup_{j=k}^{|I|} \{P_{j:I} \leq \alpha_{j:I}\})$ when the intersection hypothesis $H_I$ is true. If the Type I error probability is $\leq \alpha$, we call the local test as a level $\alpha$ test. Now, consider the family of the local tests $\{(I : \alpha_{k:I}, \ldots, \alpha_{|I|:I}) : I \subset \{1, \ldots, n\}, |I| \geq k\}$. We call this family to be symmetric if for any subsets $I$ and $J$ with $|I| = |J| \geq k$, we have $\alpha_{i:I} = \alpha_{i:J}$, for $k \leq i \leq |I|$. The notation implies that we use the same local test for testing different intersection null hypotheses $H_I$ and $H_J$ if the cardinalities of $I$ and $J$ are the same. A symmetric family of local tests is indeed characterized by a double-integer-indexed set of critical constants. Later, for simplicity, we use $(I : \alpha_{k,|I|}, \ldots, \alpha_{|I|,|I|})$ to denote the local test $(I : \alpha_{k:I}, \ldots, \alpha_{|I|:I})$ in the symmetric family. The following lemma plays an

important role in the construction of symmetric local tests.

LEMMA 3.1. (Röhmel and Streitberg, 1987; Falk, 1989; Lehmann and Romano, 2005) *Suppose $P_1, \ldots, P_t$ are p-values in the sense that $P\{P_i \leq u\} \leq u$ for all $i$ and $u$ in $(0,1)$. Let their ordered values be $P_{(1)} \leq \cdots \leq P_{(t)}$. Let $0 = \beta_0 \leq \beta_1 \leq \cdots \leq \beta_m \leq 1$ for some $m \leq t$. Then*

$$P\left\{P_{(1)} \leq \beta_1, \ldots, P_{(m)} \leq \beta_m\right\} \leq t \sum_{i=1}^{m} (\beta_i - \beta_{i-1})/i \ .$$

For any subset $I \subset \{1, \ldots, n\}$ with $|I| \geq k$, consider the local test $(I : \alpha_{k,|I|}, \cdots, \alpha_{|I|,|I|})$ of $H_I$. Let $m = t = |I|$ and $\beta_i = \alpha_{i,|I|}$ if $k \leq i \leq |I|$ and otherwise 0. Then, by Lemma 3.1, the Type I error rate of the local test is less than or equal to

$$(3.2) \qquad |I|\Big(\alpha_{k,|I|}/k + \sum_{i=k+1}^{|I|} (\alpha_{i,|I|} - \alpha_{i-1,|I|})/i\Big) \ .$$

Evidently, if the right side of (3.2) is bounded above by $\alpha$, then the local test is a level $\alpha$ test under arbitrary dependence of p-values.

Let $\alpha_k \leq \cdots \leq \alpha_n$ be given. Define

$$(3.3) \quad D_1(k) = \max_{k \leq |I| \leq n} \left\{ |I|\frac{\alpha_{n-|I|+k}}{k} + |I| \sum_{k<j\leq|I|} \frac{\alpha_{n-|I|+j} - \alpha_{n-|I|+j-1}}{j} \right\} \ .$$

Then by (6), the local test $(I : \alpha\alpha_{n-|I|+i}/D_1(k),\ k \leq i \leq |I|)$ is a level $\alpha$ tests of $H_I$. Specially, suppose $\alpha_{j,|I|}$ $(j = k, \ldots, |I|)$ is a constant. For the validity of $(3.2) \leq \alpha$, one only needs $\alpha_{k,|I|} \leq k\alpha/|I|$. Thus the local test $(I : k\alpha/|I|, \cdots, k\alpha/|I|)$ is also a level $\alpha$ test of $H_I$.

After obtaining the symmetric families of local tests, we now generalize the usual closure principle for controlling the $k$-FWER.

THEOREM 3.1. *Let $\{H_i, i = 1, \ldots, n\}$ be a finite family of hypotheses and $k$ be pre-specified with $1 \leq k \leq n$. For any $I \subset \{1, \ldots, n\}$ with $|I| \geq k$, let $P_i, i \in I$ be the p-values associated with $H_i, i \in I$ and $P_{1:I} \leq \cdots \leq P_{|I|:I}$ an ordered arrangement of these p-values. Let $H_I = \bigcap_{j \in I} H_j$. Suppose a level $\alpha$ local test defined by (3.1) is available to test $H_I$ for each $I$. Then, the (generalized) closed testing procedure, which rejects any hypothesis $H_i$ if and only if $H_I$ is rejected for all $I$ satisfying $i \in I$ and $P_i \geq P_{k:I}$, controls the k-FWER at level $\alpha$.*

PROOF. Let $I_0$ be the set of indices of true hypotheses. Assume $|I_0| \geq k$ or there is nothing to prove. Define the event

$$A = \{\text{at least } k \text{ true hypotheses are rejected}\}.$$

The occurrence of event $A$ implies that there exists $i \geq k$ such that the null hypothesis $H_{i:I_0}$ is rejected. From the description of the generalized closed testing procedure, $H_{i:I_0}$ is rejected implies that $H_{I_0}$ is rejected. Therefore,

$$k\text{-FWER} = P(A) \leq P\{H_{I_0} \text{ rejected}\} \leq \alpha.$$

□

REMARK 3.1. If $k = 1$, the generalized closed testing procedure described above is identical to the one proposed by Marcus et al. (1976). Under such closed testing procedure, $H_{(1)}, \ldots, H_{(k-1)}$ will always be rejected.

REMARK 3.2. Compared with the original closure principle, the generalized closed testing procedure involves far fewer single tests for testing significance of $H_i$ when $k \geq 2$. One of the reasons is that, when testing $H_I$, $(k-1)$ false rejections are allowed. If $P_i$ is one of the first $(k-1)$ minimal $p$-values in $P_j, j \in I$, the $p$-value $P_i$ does not play any role in testing $H_I$.

In what follows, the closed testing procedures considered are always built on symmetric families of local tests characterized by the critical constants $\alpha_{i,|I|}$.

**4. Stepwise Procedure.** In this section, we discuss how to apply the generalized closure principle enunciated in Theorem 3.1 to derive stepwise (stepup or stepdown) procedures with the $k$-FWER controlling property.

It is generally not easy to show directly that a specific stepwise procedure has the $k$-FWER controlling property. However, our strategy is first to build a closed testing procedure based on a family of level $\alpha$ local tests, and then to prove that the specific stepwise procedure is equivalent to or dominated by the closed testing procedure. We now qualify *equivalence* or *dominance* of two procedures (Liu, 1996; Grechanovsky and Hochberg, 1999). Two procedures $A$ and $B$ are called *equivalent* if they reject or accept the same *individual* hypotheses. Procedure $A$ is said to *dominate* Procedure $B$ if $A$ always rejects at least those individual hypotheses rejected by Procedure $B$. It is easy to see that, if Procedure $B$ is shown to be equivalent to or dominated by Procedure $A$, which has the $k$-FWER controlling property, then Procedure $B$ also has the $k$-FWER controlling property. We now focus on stepdown procedures.

4.1. *Stepdown Procedure.* When the $p$-values are in any dependency structure, Hommel and Hoffmann (1987) and Lehmann and Romano (2005) pro-posed a stepdown procedure with the critical values $\alpha_i$ defined by,

$$\alpha_i = \frac{k}{n-i+k}\alpha, \quad i = k, \cdots, n, \tag{4.1}$$

with the $k$-FWER controlling property, which is a generalization of Holm's procedure (Holm, 1979). In this subsection, we provide a general result (Theorem 4.1) of deriving stepdown procedures through given closed testing procedures.

THEOREM 4.1. *Suppose $\alpha_{i,|I|}, k \leq i \leq |I|$ and $k \leq |I| \leq n$ be given, $\alpha_{i,|I|}$ is increasing in $i$ and decreasing in $|I|$. The following statements are true.*
*(i) The closed testing procedure with the critical values $\alpha_{i,|I|}$ dominates the stepdown procedure with the critical values $\alpha_{k,(n-i)+k}$, $k \leq i \leq n$.*
*(ii) Furthermore, if $\alpha_{i,|I|}$, $k \leq i \leq |I|$ are constant for each given $|I|$, then these two procedures are equivalent.*
*(iii) If $(I : \alpha_{i,|I|}, k \leq i \leq |I|)$ is a level $\alpha$ local test of $H_I$ for each $I \subset \{1, \ldots, n\}$ with $|I| \geq k$, then the stepdown procedure controls the $k$-FWER at level $\alpha$.*

PROOF. (i) We first show that, for any individual hypothesis $H_{(i)}$ with index $(i)$, which corresponds to the $i$th minimal $p$-value $P_{(i)}$, if it is rejected by the stepdown procedure, it is also rejected by the closed testing procedure.

If $i < k$, $H_{(i)}$ is automatically rejected by these two procedures. We assume $i \geq k$. If $H_{(i)}$ is rejected by the stepdown procedure, then $P_{(j)} \leq \alpha_{k,(n-j)+k}$, for all $k \leq j \leq i$. Consider any subset $I \subset \{1, \ldots, n\}$ with $(i) \in I, |I| \geq k$, and $P_{(i)} \geq P_{k:I}$. Suppose $P_{k:I} = P_{(l)}$. Then

$$k \leq l \leq i \text{ and } |I| \leq k + (n-l).$$

Consequently,

$$P_{k:I} = P_{(l)} \leq \alpha_{k,(n-l)+k} \leq \alpha_{k,|I|}.$$

That is, $H_I$ is rejected. Then, from Theorem 3.1, $H_{(i)}$ is rejected by the closed testing procedure. Therefore, the stepdown procedure with the critical values $\alpha_{k,(n-i)+k}$, $k \leq i \leq n$ is dominated by the closed testing procedure.

(ii) We now show that, when $\alpha_{i,|I|}, k \leq i \leq |I|$ are constant for each given $|I|$, if $H_{(i)}$ is rejected by the closed testing procedure, it is also rejected by the stepdown procedure.

Let $I_j = \{(1), \ldots, (k-1), (j), \ldots, (n)\}$, $k \leq j \leq i$. Note that, $(i) \in I_j, |I_j| = (n-j) + k \geq k$, and $P_{k:I_j} = P_{(j)}$. If $H_{(i)}$ is rejected by the closed testing procedure, then, from Theorem 3.1, each $H_{I_j}$ will also be rejected by the corresponding level $\alpha$ local test. Note that, $\alpha_{i,|I|}, k \leq i \leq |I|$ are constant for each given $|I|$, then $H_{I_j}$ is rejected if and only if $P_{k:I_j} \leq \alpha_{k,|I_j|}$. That is, $P_{(j)} \leq \alpha_{k,(n-j)+k}, k \leq j \leq i$. Hence, from the definition of generalized stepdown procedure, $H_{(i)}$ is rejected by the stepdown procedure with the critical values $\alpha_{k,(n-i)+k}$. So, these two procedures are equivalent.

(iii) Since $(I : \alpha_{i,|I|}, k \leq i \leq |I|)$ is a level $\alpha$ local test for each $I \subset \{1, \ldots, n\}$ with $|I| \geq k$, then, from Theorem 3.1, the closed testing procedure controls $k$-FWER at level $\alpha$. Note that, the stepdown procedure is dominated by the closed testing procedure, so the stepdown procedure also controls the $k$-FWER at level $\alpha$. □

Theorem 4.1 shows that under certain conditions, a generalized closed testing procedure can be reduced to a corresponding stepdown procedure. As an application, we point out that Hommel and Hoffmann (1987) and Lehmann and Romano (2005)'s stepdown procedure can be derived by the general result by choosing appropriate level $\alpha$ local tests. As illustrated in Section 3, the local test $(I : \frac{k}{|I|}\alpha, \ldots, \frac{k}{|I|}\alpha)$ is a level $\alpha$ test of $H_I$ for each $I$. Note that $\alpha_{i,|I|} = \frac{k}{|I|}\alpha$, $k \leq i \leq |I|$ are constant for each given $|I|$. Thus, by Theorem 4.1, the closed testing procedure based on these local tests is equivalent to the stepdown procedure with the critical values defined by (4.1).

We now focus on the converse of Theorem 4.1. In Theorem 4.2, we show that any stepdown procedure is equivalent to some closed testing procedure.

THEOREM 4.2. *Let $\alpha_k \leq \cdots \leq \alpha_n$ be given. Then, the stepdown procedure with the critical values $\alpha_i$, $k \leq i \leq n$ is equivalent to a closed testing procedure with the critical values $\beta_{i,|I|} = \alpha_{n-|I|+k}$, $k \leq i \leq |I|$ and $k \leq |I| \leq n$.*

PROOF. Note that, $\beta_{i,|I|} = \alpha_{n-|I|+k}$ is increasing in $i$ and decreasing in $|I|$, and for each given $|I|$, $\beta_{i,|I|}$ are constant. Then, from Theorem 4.1, the closed testing procedure with critical values $\beta_{i,|I|} = \alpha_{n-|I|+k}$ is equivalent to the stepdown procedure with the critical values $\beta_{k,n-i+k} = \alpha_i$. $\square$

Theorem 4.2 shows that the stepdown procedures can be viewed as a specific class of the closed testing procedures.

4.2. *Stepup Procedure.* When the $p$-values are in any dependency structure, Romano and Shaikh (2006b) obtained a stepup procedure with the critical values $\alpha'_i = \alpha\alpha_i/D_1(k)$, where $\alpha_k \leq \cdots \leq \alpha_n$ are any non-decreasing constants and $D_1(k)$ is defined as in (3.3). In this subsection, we give a general result (Theorem 4.3) of deriving stepup procedures through given closed testing procedures.

THEOREM 4.3. *Suppose $\alpha_{i,|I|}$, $k \leq i \leq |I|$ and $k \leq |I| \leq n$ be given, $\alpha_{i,|I|}$ is increasing in $i$ and decreasing in $|I|$, and $\alpha_{l,(n-i)+l}$ is increasing in $l$ for each given $i \geq k$. The following statements are true.*

*(i) The closed testing procedure with the critical values $\alpha_{i,|I|}$ dominates the stepup procedure with the critical values $\alpha_{k,(n-i)+k}$.*

*(ii) Furthermore, if $\alpha_{l,(n-i)+l}$, $k \leq l \leq n$ are constant for each given $i \geq k$, then these two procedures are equivalent.*

*(iii) If $(I : \alpha_{i,|I|}, k \leq i \leq |I|)$ is a level $\alpha$ local test of $H_I$ for each $I \subset \{1,\ldots,n\}$ with $|I| \geq k$, then the stepup procedure controls the k-FWER at level $\alpha$.*

PROOF. (i) We show that, for any individual hypothesis $H_{(i)}$, if it is rejected by the stepup procedure, it is also rejected by the closed testing procedure.

If $i < k$, $H_{(i)}$ is automatically rejected by these two procedures. So, we assume $i \geq k$. If $H_{(i)}$ is rejected by the stepup procedure, then there exists $j \geq i$ satisfying $P_{(j)} \leq \alpha_{k,(n-j)+k}$. Consider any subset $I \subset \{1,\ldots,n\}$ with $(i) \in I$, $|I| \geq k$, and $P_{(i)} \geq P_{k:I}$. Let $l = \max\{i' : P_{i':I} \leq P_{(j)}\}$. Since $j \geq i$ and $(i) \in I$, $l \geq k$. If $l < |I|$, then $P_{l+1:I} > P_{(j)}$, so $|I| \leq (n-j)+l$. Evidently, when $l = |I|$, it follows that $|I| \leq (n-j)+l$, too. Hence,

$$P_{l:I} \leq P_{(j)} \leq \alpha_{k,(n-j)+k} \leq \alpha_{l,(n-j)+l} \leq \alpha_{l,|I|}.$$

The third inequality in the chain above follows from the assumption that $\alpha_{l,(n-i)+l}$ is increasing in $l$, and the last inequality follows from the inequality $|I| \leq (n-j)+l$. Consequently, $H_I$ is rejected. Then, by Theorem 3.1, $H_{(i)}$ is rejected by the closed testing procedure. Hence, the stepup procedure with the critical values $\alpha_{k,(n-i)+k}$ is dominated by the closed testing procedure.

(ii) We show that when $\alpha_{l,(n-j)+l}$, $k \leq l \leq n$ are constant for each given $j \geq k$, if $H_{(i)}$ is rejected by the closed testing procedure, it is also rejected by the stepup procedure.

Let $I = \{(1), \ldots, (k-1), (i), \ldots, (n)\}$, and note that, $(i) \in I, |I| = (n-i) + k \geq k$. If $H_{(i)}$ is rejected by the closed testing procedure, then by Theorem 3.1, $H_I$ will be rejected by the corresponding level $\alpha$ local test. That is, there exists $j \geq i$ satisfying $P_{j:I} \leq \alpha_{j,|I|}$. Since $P_{j:I} = P_{(i+j-k)}$, $\alpha_{j,|I|} = \alpha_{j,n-i+k}$, and $\alpha_{l,(n-i')+l}$, $k \leq l \leq n$ are constant for each given $i' \geq k$, $\alpha_{j,n-i+k} = \alpha_{k,n-i-j+2k}$. Hence, $P_{(i+j-k)} \leq \alpha_{k,n-(i+j-k)+k}$. From the definition of generalized stepup procedure, $H_{(i)}$ is rejected by the stepup procedure. Consequently, these two procedures are equivalent.

(iii) Since $(I : \alpha_{i,|I|}, k \leq i \leq |I|)$ is a level $\alpha$ local test of $H_I$ for each $I \subset \{1, \ldots, n\}$ with $|I| \geq k$, by Theorem 3.1, the closed testing procedure controls the $k$-FWER at level $\alpha$. Note that the stepup procedure is dominated by the closed testing procedure, so the stepup procedure also controls the $k$-FWER at level $\alpha$. □

Theorem 4.3 shows that under certain conditions, a generalized closed testing procedure can be reduced to a corresponding stepup procedure. As an application, we show that Romano and Shaikh (2006)'s stepup procedure can be derived by the general result by choosing appropriate level $\alpha$ local tests. Let $\alpha_k \leq \cdots \leq \alpha_n$ be given and define $\beta_{i,|I|} = \alpha_{n-|I|+i}$. As illustrated in Section 3, $(I : \alpha\alpha_{n-|I|+k}/D_1(k), \cdots, \alpha\alpha_n/D_1(k))$ is a level $\alpha$ local test of $H_I$ for each $I$. Note that $\beta_{l,(n-i)+l} = \alpha_i$, $k \leq l \leq n$ are constant for each given $i$. Then, based on these local tests, Romano and Shaikh (2006)'s stepup procedure is derived by Theorem 4.3.

We now focus on the converse of Theorem 4.3. In Theorem 4.4, we show that any stepup procedure is equivalent to some closed testing procedure.

THEOREM 4.4. *Let $\alpha_k \leq \cdots \leq \alpha_n$ be given. Then, the stepup procedure with the critical values $\alpha_i$, $k \leq i \leq n$ is equivalent to a closed testing procedure with the critical values $\beta_{i,|I|} = \alpha_{n-|I|+i}$, $k \leq i \leq |I|$ and $k \leq |I| \leq n$.*

PROOF. Note that, $\beta_{i,|I|} = \alpha_{n-|I|+i}$ is increasing in $i$, decreasing in $|I|$, and $\beta_{l,n-i+l} = \alpha_i$ are constant for each given $i \geq k$. Then, from Theorem 4.3, the closed testing procedure with critical values $\beta_{i,|I|} = \alpha_{n-|I|+k}$ is equivalent to the stepup procedure with the critical values $\beta_{k,n-i+k} = \alpha_i$, $k \leq i \leq n$. □

Theorem 4.4 shows that the stepup procedures can also be viewed as a specific class of the closed testing procedures.

REMARK 4.1. When the $p$-values are positively dependent in the sense of being multivariate totally positive of order two (MTP$_2$) (Karlin and Rinott, 1980), based on the generalized Simes' test of Sarkar (2008a), we can easily derive the stepdown and stepup procedures in Sarkar (2008a) by Theorems 4.1 and 4.3.

**5. Generalized Hommel Procedure.** In the context of the FWER, Hommel (1988) developed a well-known sequential procedure based on Simes (1986)' test, which is more powerful than Hochberg (1988)'s stepup procedure (Hommel, 1989). Hommel's procedure is described as follows: compute $j = \max\{i \in \{1, \ldots, n\} : P_{(n-i+l)} > \frac{l\alpha}{i}, \text{ for } l = 1, \ldots, i\}$. If the maximum does not exist, reject all $H_i$ $(i = 1, \ldots, n)$; otherwise reject all $H_i$ with $P_i \leq \alpha/j, (i = 1, \ldots, n)$. In this section, we generalize Hommel's procedure in two directions. In one direction, we move from the critical values of Simes' test to any double-indexed critical values satisfying certain properties, and in another, we move from the control of the FWER to that of the $k$-FWER.

For convenience of discussion, a formal definition of a generalized Hommel procedure is first given as follows.

DEFINITION 5.1. *Let $\alpha_{i,|I|}, k \leq i \leq |I|$ and $k \leq |I| \leq n$ be given. Suppose $\alpha_{i,|I|}$ is increasing in $i$ and decreasing in $|I|$. A generalized Hommel procedure with the critical values $\alpha_{i,|I|}$ is defined as follows: compute $j = \max\{i \in \{k, \ldots, n\} : P_{(n-i+l)} > \alpha_{l,i}, \text{ for } l = k, \ldots, i\}$. If the maximum does not exist, reject all $H_{(i)}$ $(i = k, \ldots, n)$, otherwise reject all $H_{(i)}$ with $P_{(i)} \leq \alpha_{k,j}$ $(i = k, \ldots, n)$. In any case, the first $(k-1)$ most significant hypotheses $H_{(1)}, \ldots, H_{(k-1)}$ are automatically rejected.*

For example, for the original Hommel procedure, $k = 1$, and $\alpha_{i,|I|}$ is taken to be $\alpha_{i,|I|} = \frac{i}{|I|}\alpha$. The following is the main result of this section.

THEOREM 5.1. *Let $\alpha_{i,|I|}, k \leq i \leq |I|$ and $k \leq |I| \leq n$ be given, and suppose $\alpha_{i,|I|}$ is increasing in $i$ and decreasing in $|I|$. Then the following are true.*
*(i) The closed testing procedure with the critical values $\alpha_{i,|I|}$ is equivalent to the generalized Hommel procedure with the same critical values.*
*(ii) If $(I : \alpha_{i,|I|}, \ k \leq i \leq |I|)$ is a level $\alpha$ local test of $H_I$ for each $I \subset \{1, \cdots, n\}$ with $|I| \geq k$, then the generalized Hommel procedure controls the $k$-FWER at level $\alpha$.*

PROOF. (i) If $i < k$, $H_{(i)}$ is automatically rejected by these two procedures. So, we assume $i \geq k$. Let $j = \max\{i' \in \{k, \ldots, n\} : P_{(n-i'+l)} > \alpha_{l,i'}, \text{for } l = k, \ldots, i'\}$. Suppose $j$ does not exist. By Definition 5.1, the generalized Hommel procedure rejects all $H_i$. We now show that the closed testing procedure rejects all $H_i$. Non-existence of $j$ implies that for any $i$ with $k \leq i \leq n$, there exists $l$ with $k \leq l \leq i$ satisfying $P_{(n-i+l)} \leq \alpha_{l,i}$. Consider any subset $I \subset \{1, \ldots, n\}$ with $|I| = i$. Suppose $P_{l:I} = P_{(l')}$, then $|I| \leq (n - l') + l$. That is, $l' \leq n - i + l$. So, $P_{l:I} \leq P_{(n-i+l)} \leq \alpha_{l,|I|}$. Hence, $H_I$ is rejected. By Theorem 3.1, all hypotheses $H_{(1)}, \ldots, H_{(n)}$ are rejected by the closed testing procedure.

We now consider the case that $j$ exists. First, we show that, for any individual hypothesis $H_{(i)}$, if it is rejected by the generalized Hommel procedure, it is also rejected by the closed testing procedure.

If $H_{(i)}$ is rejected by the generalized Hommel procedure, then $P_{(i)} \leq \alpha_{k,j}$. Consider any subset $I \subset \{1, \ldots, n\}$ with $(i) \in I, |I| \geq k$, and $P_{(i)} \geq P_{k:I}$. If $|I| \leq j$, then

$$P_{k:I} \leq P_{(i)} \leq \alpha_{k,j} \leq \alpha_{k,|I|}.$$

Consequently, $H_I$ is rejected. If $|I| > j$, then from the definition of $j$, there exists $l$ satisfying $k \leq l \leq |I|$ such that $P_{(n-|I|+l)} \leq \alpha_{l,|I|}$. Let $i_0 = \max\{i' \geq k : P_{i':I} \leq P_{(n-|I|+l)}\}$. Note that $(|I| - k + 1) + (n - |I| + l) = n + l - k + 1 > n$, which implies that the maximum $i_0$ exists. From the definition of $i_0$, we have $P_{i_0+1:I} > P_{(n-|I|+l)}$, so $|I| \leq n - (n - |I| + l) + i_0$. That is, $l \leq i_0$. Therefore,

$$P_{i_0:I} \leq P_{(n-|I|+l)} \leq \alpha_{l,|I|} \leq \alpha_{i_0,|I|}.$$

Hence, $H_I$ is rejected. Consequently, $H_{(i)}$ is rejected by the closed testing procedure. Therefore, the generalized Hommel procedure with the critical values $\alpha_{i,|I|}$ is dominated by the closed testing procedure.

Next, we show that, if $H_{(i)}$ is rejected by the closed testing procedure, it is also rejected by the generalized Hommel procedure.

Let $I = \{(1), \ldots, (k-1), (n-j+k), (n-j+k+1), \ldots, (n)\}$. Note that $|I| = j \geq k$, and $P_{l:I} = P_{(n-j+l)}, k \leq l \leq j$. From the definition of $j$, we have $P_{(n-j+l)} > \alpha_{l,j}$ for all $k \leq l \leq j$. Thus, $H_I$ will not be rejected by the local test $(I : \alpha_{i,j}, k \leq i \leq j)$. If $H_{(i)}$ is rejected by the closed testing procedure, then from Definition 5.1, $(i) \notin I$. That is, $k \leq i < n - j + k$. Let $I' = \{(1), \ldots, (k-1), (i), (n-j+k+1), \ldots, (n)\}$. Note that $(i) \in I'$, $|I'| = j \geq k$ and $P_{k:I'} = P_{(i)}$. Then $H_{I'}$ is rejected following from Theorem 3.1. It now follows that $P_{(i)} \leq \alpha_{k,j}$ and $H_{(i)}$ is rejected by the generalized Hommel procedure. Hence, these two procedures are equivalent.

(ii) Since $(I : \alpha_{i,|I|}, k \leq i \leq |I|)$ is a level $\alpha$ local test for each $I \subset \{1, \ldots, n\}$ with $|I| \geq k$, then, from Theorem 3.1, the closed testing procedure controls the $k$-FWER at level $\alpha$. Since the generalized Hommel procedure is equivalent to the closed testing procedure, the generalized Hommel procedure also controls the $k$-FWER at level $\alpha$. □

Theorem 5.1 shown that each generalized closed testing procedure is equivalent to a generalized Hommel procedure with the same critical values. Since generalized Hommel procedures are much simpler than closed testing procedures, this result implies that we find a new shortcut, which can reduce generalized closed testing procedures to computationally simple procedures in applications.

REMARK 5.1. Professor Hommel brought to our attention, while commenting on an earlier draft of the paper, that a version similar to Theorem 5.1 for $k = 1$ appears in the paper by Bernhard et al. (2004). Their result is stated without proof.

Finally, we provide an intuitive interpretation for the generalized Hommel procedures. For each generalized Hommel procedure with the critical values $\alpha_{i,|I|}$, the value $\hat{j} = \max\{i \in \{k, \ldots, n\} : P_{(n-i+l)} > \alpha_{l,i}, \text{ for } l = k, \ldots, i\}$ can be viewed as an estimate of the number of true null hypotheses. For example, suppose $\alpha_{i,|I|} = \frac{k}{|I|}\alpha$. Let $\hat{j} = \max\left\{i \in \{k, \cdots, n\} : P_{(n-i+k)} > \frac{k}{i}\alpha\right\}$, then $\hat{j}$ can be expressed as $n - \widetilde{j} + k$, where $\widetilde{j} = \min\left\{j \in \{k, \ldots, n\} : P_{(j)} > \frac{k}{n-j+k}\alpha\right\}$. Note that $\widetilde{j} - 1$ is the number of rejected hypotheses by using the stepdown procedure of Hommel and Hoffmann (1987) and Lehmann and Romano (2005), and $(k - 1)$ is the maximal number of false rejections that one is willing to tolerate in this procedure. Thus we use $n - (\widetilde{j} - 1) + (k - 1) = n - \widetilde{j} + k = \hat{j}$ as an estimate of the number of true null hypotheses in the family of hypotheses $H_1, \ldots, H_n$. Hence, each generalized Hommel procedure with the critical values $\alpha_{i,|I|}$ can be interpreted as a two-stage procedure, in which, first estimate the number of true null hypotheses by using $\hat{j}$, and then based on the estimate $\hat{j}$, establish a single-step procedure with the critical value $\alpha_{k,\hat{j}}$.

**6. Concluding Remarks.** The original closure principle was formulated by Marcus et al. (1976) in the context of the FWER and has since been a powerful tool for deriving multiple testing procedures controlling the FWER. In fact, almost all FWER controlling procedures are either derived using this principle or can be rewritten as associated closed testing procedures. The only disadvantage is that closed testing procedures are computationally complex.

In this paper, we have generalized the closure principle for the $k$-FWER. In the same vein as the usual closure principle, the value of generalized closure principle is that, instead of simultaneously testing multiple hypotheses, one only needs to do a number of single tests of intersection null hypotheses. The construction of single tests based on marginal $p$-values is relatively easy. The reason is that, as all $p$-values associated with true null hypotheses are marginally stochastically dominated by uniformly distributed random variables, some powerful probability inequalities on uniformly distributed random variables are available such as Bonferroni inequality, Simes' inequality and generalized Simes' inequality, which are useful in building single tests. See Sarkar (2008a).

The generalized closed testing procedures are also computationally complex. We have discussed how to reduce generalized closed testing procedures to simple stepwise procedures, and showed that, under certain conditions, a generalized closed testing procedure can be formulated as a stepwise procedure. We have generalized the well-known Hommel procedure, and shown that each generalized closed testing procedure is equivalent to a simple generalized Hommel procedure with the same critical values.

We need to point out that, in this paper, we only discussed how to derive the $k$-FWER controlling procedures based on marginal $p$-values using the generalized closure principle. A future research is to use this principle for developing some resampling-based methods such as those described in Dudoit et al. (2004), van der Laan et al. (2005), and Romano and Wolf (2007) that take into account the dependence structure of the underlying test statistics.

We also need to note that the generalized closure principle is derived for the families of non-hierarchical null hypotheses. In many cases, we need to test families of hierarchical null hypotheses simultaneously. An interesting future research is to modify the generalized closure principle for families of hierarchical null hypotheses extending the work of Hommel (1986) on the FWER to the $k$-FWER.

**Acknowledgements.** The first author wishes to thank Gerhard Hommel and Joseph Romano for helpful discussions.**References.**

Biostatistics Branch
National Institute of Environmental Health Sciences
Research Triangle Park, NC 27709, U.S.A.
Department of Environmental Health
University of Cincinnati
Cincinnati, OH 45267, USA
E-mail: wenge.guo@gmail.com
raomb@ucmail.uc.edu